\documentclass[12pt]{article}
\usepackage[latin1]{inputenc}
\usepackage{amssymb}
\usepackage{amsmath}
\usepackage{amsthm}

\input pictex.tex
\newcommand{\esp}{\hspace{0.1cm}}

\newcommand{\R}{\mathbb{R}}

\newcommand{\N}{\mathbb N}
\newcommand{\Z}{\mathbb Z}

\theoremstyle{definition}

\newtheorem{thm}{Theorem}[section]

\newtheorem{prop}[thm]{Proposition}

\newtheorem{lem}[thm]{Lemma}

\newtheorem{rem}[thm]{Remark}

\newtheorem{ex}[thm]{Example}

\newcommand{\vsp}{\vspace{0.1cm}}

\date{}
\author{Crist\'obal Rivas}

\begin{document}

\title{On spaces of Conradian group orderings}
\maketitle

\vspace{-0.5cm}

\begin{abstract} \noindent We classify $C$-orderable groups admitting only finitely many
$C$-orderings. We show that if a $C$-orderable group has infinitely many $C$-orderings,
then it actually has uncountably many $C$-orderings, and none of these is isolated in
the space of $C$-orderings. As a relevant example, we carefully study the case of
Baumslag-Solitar's group $B(1,2)$. We show that $B(1,2)$ has four $C$-orderings,
each of which is bi-invariant, but its space of left-orderings is homeomorphic to
the Cantor set.
\end{abstract}

\vspace{0.4cm}

\noindent{\large \bf Introduction}

\vspace{0.532cm}

One of the starting points of the theory of orderable groups is
\cite{holder}, where O. H\"older proved that any Archimedean Abelian
ordered group is ordered isomorphic to a subgroup of the additive
group of real numbers with the standard ordering. In his seminal
work \cite{conrad}, P.~Conrad obtained a condition on left-ordered
groups which is equivalent to the fact that the conclusion of
H\"older's theorem holds `locally' (see (4) below). Since then,
these so-called $C$-orderings (or Conradian orderings) have played a
fundamental role in the theory of left-orderable groups. (See for
instance \cite{witte,rr}.) Recall that a left-invariant (total)
ordering $\, \preceq \,$ on a group $G$ is said to be
\textit{Conradian} if the following four equivalent properties hold
(this equivalence will be referred to as the Conrad Theorem, see
\cite{botto,conrad,kopytov}):

\vsp\vsp

\noindent (1) For all $f \succ id$ and $g \succ id \; $ (for all
{\em positive} $f,g$, for short), we have $f g^n \succ g$ for some
$n \in \mathbb{N}$.

\vsp\vsp

\noindent (2) If $1 \prec g  \prec f$, then $g^{-1} f^n g \succ f$
for some $n \in \mathbb{N} = \{1,2, \ldots \}$.

\vsp\vsp

\noindent (3) For all positive $g \in G$, the set $S_g = \{f \in G \mid
\esp f^n \prec g, \mbox{ for all } n \in \mathbb{Z} \}$ is a convex subgroup.
(By definition, a subset $S \subset G$ is {\em convex} if whenever \esp
$f_1 \prec h \prec f_2$ \esp for some $f_1, f_2$ in $S$, we have $h \in S$.)

\vsp\vsp

\noindent (4) Given $g \in G$, we denote the maximal (resp. minimal)
convex subgroup which does not contain (resp. contains) $g$ by
$G_g$ (resp. $G^g$). For every $g$, we have that $G_g$ is normal
in $G^g$, and there exists a non-decreasing group homomorphism
(to be referred to as the {\em Conrad homomorphism}) \esp
$\tau_{\preceq}^{g} \!: G^g \rightarrow \mathbb{R}$ \esp whose kernel
coincides with $G_g$. Moreover, this homomorphism is unique up to
multiplication by a positive real number.

\vsp\vsp\vsp

Recently, two new approaches to this property have been proposed by
A. Navas. On the one hand, as it was noticed in \cite{leslie,navas},
in (1) and (2) above one may actually take $n \!=\! 2$. The
topological counterpart of this is the fact that the space of
$C$-orderings is compact when it is endowed with a natural topology
(see \S \ref{pre-1}). This leads, for instance, to a new and short
proof of the fact that locally indicable groups are $C$-orderable.
(Note that the converse follows from (4).) On the other hand, the
dynamical characterization of the Conradian property of
\cite{navas,crossings} leads to applications in the study of the
topology of space of group orderings, and to general `level
structure' theorems for left-ordered groups. In this work, this
dynamical point of view will be crucial.

\vsp

Following the first direction above, we focus on the
structure of the space of group $C$-orderings. In particular,
we provide complete answers to questions in \cite[Question
3.9]{navas} and \cite[\S 1.3]{thompson}.

\vsp

It is known that the space of left-orderings of a group is either finite
or uncountable \cite{linnell,crossings}. Although this is no longer true
for bi-orderings \cite{butts}, our first main result shows that this
dichotomy persists for $C$-orderings.

\vspace{0.53cm}

\noindent{\bf Theorem A.} {\em Let $G$ be a $C$-orderable group. If
$G$ admits infinitely many $C$-orderings, then it has uncountably
many $C$-orderings. Moreover, none of these is isolated in the space of
$C$-orderings.}

\vspace{0.53cm}

For the second claim of Theorem A, the space of group left-orderings
is endowed with the projective topology induced from the discrete
one on finite sets, and the subset of $C$-orderings is endowed with
the subspace one (see \S \ref{pre-1} for more details). The space of
left-orderings is (Hausdorff, totally disconnected and) compact, and
the subset of $C$-orderings is closed therein. In particular, this
implies that the second claim of the statement is stronger than
the first.

\vsp

Our second result concerns groups admitting only finitely many
$C$-orderings, and may be considered as an analogue of Tararin's
classification of left-orderable groups admitting finitely many
left-orderings \cite[Theorem 5.2.1]{kopytov}. For the statement,
recall that a series
$$\{ id \} = G_0 \lhd G_{1} \lhd \ldots \lhd G_{n-1} \lhd G_n = G$$
is said to be {\em rational} if it is subnormal ({\em i.e.,} each $G_{i}$ is normal
in $G_{i+1}$) and each quotient $G_{i+1} / G_{i}$ is torsion-free rank-1 Abelian.

\vspace{0.53cm}

\noindent{\bf Theorem B.} {\em Let $G$ be a $C$-orderable group. If
$G$ admits finitely many $C$-orderings, then $G$ admits a unique
(hence normal) rational series. In this series, no quotient
$G_{i+2}/G_{i}$ is Abelian. Conversely, if $G$ is a group admitting
a  normal rational series
$$\{ id \} = G_0 \lhd G_{1} \lhd \ldots \lhd G_{n-1} \lhd G_n = G$$
so that no quotient $G_{i+2} / G_{i}$ is Abelian, then the number
of $C$-orderings on $G$ equals $2^n$.}

\vspace{0.53cm}

The proof of Theorem B consists in a non-trivial modification of
Tararin's arguments. (Note that the statement of Tararin's theorem
is the same as that of Theorem B though changing `$C$-orderings' by
`left-orderings', and the condition `$G_{i+2} / G_{i}$ non Abelian'
by `$G_{i+2} / G_{i}$ non bi-orderable'.) First, as observed
in \cite{thompson}, if a $C$-orderable group admits finitely
many $C$-orderings, then it must be solvable. Now the fact
that each quotient $G_{i+1} / G_{i}$ has rank 1 and no quotient
$G_{i+2} / G_{i}$ is Abelian is a consequence of the fact that the space
of orderings of higher-rank torsion-free Abelian groups are uncountable
(see for example \cite{clay, sikora}). Finally, we use an extra argument
involving the Conrad homomorphism to show that, when the hypotheses are
fulfilled, there are only finitely many $C$-orderings.

\vsp

Let us point out that Baumslag-Solitar's group \esp $B(1,2) =
\langle a, b \mid b a b^{-1} = a^2 \rangle$ \esp satisfies the
conditions of Theorem B. Therefore, its space of $C$-orderings
is finite. Actually, this space consists of four $C$-orderings,
each of which is bi-invariant (see Proposition \ref{solo-cuatro}).
This example was the starting point of this work,
and we provide a direct short argument for this particular
case in \S \ref{ejem}. We point out, however, that the space
of left-orderings of $B(1,2)$ is uncountable. (Actually, it is
homeomorphic to the Cantor set.) In addition, in \S \ref{ejem} we
give a complete description of all left-orderings of $B(1,2)$
which extends to many other left-orderable metabelian groups
(see Theorem \ref{laprop}).

\vsp

This work corroborates a general principle concerning $C$-orderings. On
the one hand, these are sufficiently rigid in that they allow deducing
structure theorems for the underlying group ({\em e.g.}, local
indicability). However, they are still sufficiently malleable in
that, starting with a $C$-ordering on a group, one may create very
many $C$-orderings, which turn out to be different from the original
one (see Example \ref{exflip}) with the only exception of the pathological
cases described in Theorem B.


\section{Preliminaries}

\subsection{Spaces of group orderings}
\label{pre-1}

\hspace{0.45cm} Given a left-orderable group $G$ (of arbitrary
cardinality), we denote the set of all left-orderings on $G$ by
$\mathcal{LO} (G)$. This set has a natural topology: a basis
of neighborhoods of $\,\preceq \,$ in $\, \mathcal{LO} (G)$ is the family of
the sets $\, U_{g_1,\ldots,g_k}$ of all left-orderings $\,\preceq'\,$ on
$G$ which coincide with $\, \preceq \,$ on $\{g_1,\ldots,g_k\}$, where
$\{g_1,\ldots,g_k\}$ runs over all finite subsets of $G$. Another
basis is given by the sets $\, V_{f_1,\ldots,f_k} \,$ of all left-orderings
$\, \preceq' \,$ on $G$ such that all the $\,f_i\,$ are $\,\preceq'$-positive, where
$\{f_1,\ldots,f_k\}$ runs over all finite subsets of $\,\preceq$-positive
elements of $G$. Endowed with this topology, $\,\mathcal{LO}(G)$ is Hausdorff and totally
disconnected, and by (an easy application of) the Tychonov Theorem,
it is compact (see for instance \cite[\S 2.1]{navas}).
The (perhaps empty) subspaces $\, \mathcal{BO}(G) \,$ and
$\, \mathcal{CO}(G) \,$ of bi-orderings and $C$-orderings on $G$
are, respectively, closed inside $\mathcal{LO}(G)$, hence
compact.

\vsp

If $G$ is countable, then this topology is metrizable: given an
exhaustion $G_0 \subset G_1 \subset \ldots$ of $G$ by finite sets,
for different $\, \preceq \,$ and $\,\preceq'\,, \,$ we may define
$\,dist(\preceq,\preceq') = 1 / 2^n$, where $n$ is the first integer such
that $\, \preceq \,$ and $\, \preceq'\,$ do not coincide on $G_n$. If $G$ is
finitely generated, we may take $G_n$ as the ball of radius $n$
with respect to a fixed finite system of generators.

\begin{ex} In the case of Conradian orderings, there is a natural
way to generate new $C$-orderings starting with a given one. This
procedure is useful for approximating a given $C$-ordering if the
series of convex subgroups is long enough (see \S \ref{B}).
Let $\, \preceq \,$ be a $C$-ordering, and let
$$\{id\} = G^{id} \subset \ldots \subset G_g \lhd G^g \subset\ldots \subset G$$
be the (perhaps infinite) series of $\, \preceq$-convex subgroups.
Taking any $g \in G \setminus \{id\}$, we may obtain a different
$C$-ordering $\, \preceq_g \,$ by `flipping' the ordering on the quotient
$\, G^g/G_g$. More precisely, given $f\in G$, we define \esp $f\succ_g
id$ \esp if one of the following (mutually excluding) conditions
holds:

\vsp

\noindent -- \esp  $f\succ id$ and $f  \not \in G^g$,

\vsp

\noindent -- \esp $f\succ id$ and $f\in G_g$,

\vsp

\noindent -- \esp $f \prec id $ and $f\in G^g\setminus G_g$.

\vsp\vsp

Clearly, this is a total ordering. To see that it is
left-invariant, we need to check that the product of any two $\,\preceq_g$-positive
elements $h_1,h_2\,$ is still $\,\preceq_g$-positive. This is obvious if $h_1 = h_2$. Now
if $\, 1 \prec_g h_1 \prec_g h_2 \,, \,$ then it is easy to check that both $\, h_1h_2 \,$ and $\, h_2h_1 \,$ belong to $\, G^{h_2} \setminus G_{h_2}. \,$ Therefore, the $\,\preceq$-\textit{signs} of $\,h_1h_2 \,$ and $\,h_2h_1 \,$ are the same as that of $h_1$, which implies that $\, h_1 h_2 \,$ and $\, h_2 h_1 \,$ are $\, \preceq_g$-positive.

\vsp

Finally, to see that $\,\preceq_g \,$ is Conradian, it suffices to show
that $id \prec_g h_1 \preceq_g h_2$ implies $\, h_1^{-1} h_2 h_1
\succ_g id \,$ and $\, h_2^{-1} h_1 h_2^2 \succ_g id\,. \,$ The first
inequality follows from $h_1^{-1} h_2 \succeq_g id$ and $h_1
\succ_g id$ just using the fact that the product of two positive
elements is still positive. For the second inequality, note that
$h_1$ and $h_2$ commute modulo $G_{h_2}$. Therefore, $h_2^{-1} h_1
h_2^2 G_{h_2} = h_2 G_{h_2}$, which implies that $h_2^{-1} h_1 h_2^2
\succ_g id$. \label{exflip}
\end{ex}

For applications of the technique of the preceding example to the
problem of approximation of group orderings, see \cite{zenkov}.


\subsection{From ordered representations to Conradian orderings}

\hspace{0.45cm} We begin by recalling an old theorem due to P. Cohn,
M. Zaitseva, and P. Conrad (see \cite[Theorem 3.4.1]{kopytov}):

\vsp

\begin{thm}\label{Cohn} {\em A group $G$ is left-orderable if and only if it
embeds in the group of (order-preserving) automorphisms of a totally
ordered set.}
\end{thm}

\vsp

Both implications of this theorem are easy. In one direction,
note that a left-ordered group acts on itself by order preserving
automorphisms, namely left translations. Conversely, to create a
left-ordering on a group $G$ of automorphisms of a totally ordered
set $(\Omega,\leq)$, we construct the what is called
\textit{induced ordering} from the action as follows. Fix a
well-order $\, \leq^* \,$ on the elements of $\,\Omega\,, \,$ and, for every $f\in
G\,, \,$ let $\, w_f = \min_{\leq^*} \{w \in \Omega \mid f(w)\not= w\}$.
Then we define an ordering $\, \preceq \,$ on $G$ by letting $f\succ id\,$
if and only if $\;f(w_f)>w_f$. It is not hard to check that this
order relation is a (total) left-ordering on $G$.

\vsp

In what follows, we will need an important definition which was
introduced in \cite{crossings}. Let $G$ be a group acting
by order preserving bijections on a totally ordered space
$(\Omega,\leq)$. A {{\em crossing}} for the action of $G$
on $\Omega$ is a 5-uple $\, (f,g,u,v,w) \,$ where $f,g$
(resp. $u,v,w$) belong to $G$ (resp. $\Omega$) and
satisfy:

\vsp

\noindent $i)$ \esp $u < w < v$.

\vsp

\noindent $ii)$ \esp For every $\,n \in \mathbb{N} , \,$
we have $\,g^n u < v $ \esp and \esp $f^n v > u \,.$

\vsp

\noindent $iii)$ \esp There exist $M,N$ in $\mathbb{N}$ such that \esp
$f^N v < w < g^M u$.

\vsp \vsp

The reason why this definition is so important is because
it actually characterizes the $C$-orderings, as is shown
in \cite[Theorem 1.4]{crossings}. We quote the theorem below.

\vsp

\begin{thm} {\em A left-ordering $\, \preceq \,$ on $G$ is Conradian if and only
if the action of $G$ by left translations on itself admits no
crossing (when taking $\,(\Omega,\leq) = (G,\preceq)$).}
\label{thm-navas}
\end{thm}

\vsp

The following crucial lemma is essentially proved in \cite{navas} in
the case of countable groups, but the proof therein rests upon very
specific issues about the so-called {\em dynamical realization}
of an ordered group. Here we give a general algebraic proof.

\vsp

\begin{lem} \label{pullback} {\em If a faithful action of a group $G$
by automorphisms of an ordered set $\,\Omega \,$ has no crossing, then
any induced ordering on $G$ is Conradian.}
\label{generalization}
\end{lem}

\noindent{\em Proof.} Suppose that the ordering $\, \preceq \,$ on $G$
induced from some well-order $\, \leq^* \,$ on $\Omega$ is not Conradian.
Then there are $\,\preceq$-positive elements $f,g$ in $G$ such that
$fg^{n} \prec g$, for every $n \in \mathbb{N}$. This easily implies
$f \prec g$. Let $\bar{w} =
\min_{\leq^*} \{ w_f,\, w_g \}$. We claim that $(fg,fg^2, \bar{w},
g(\bar{w}), fg^2(\bar{w}))$ is a crossing (see Figure 1). Indeed,
the inequalities $id \prec f \prec g$ imply that $\bar{w}=w_g\leq^*
w_f$ and $g(\bar{w})>\bar{w}$. Moreover $f(\bar{w}) \geq \bar{w} \,, \,$
which together with $fg^{n} \prec g \,$ yield $\bar{w}< fg^2(\bar{w})<
g(\bar{w})$, hence condition~$i)$ of the definition of crossing is
satisfied. Note that the preceding argument actually shows that
$fg^n(\bar{w})<g(\bar{w})$, for all $n \in \mathbb{N}\,. \,$ Thus
$fg^2fg^2(\bar{w})< fg^3(\bar{w})< g(\bar{w})$. A straightforward
induction argument shows that $(fg^2)^n(\bar{w})< g(\bar{w})$, for
all $n \in \mathbb{N}$, which proves the first part of
condition~$ii)$. For the second part, from $g(\bar{w})>\bar{w}\,$ and
$\,f(\bar{w})\geq\bar{w} \,$ we conclude that $\, \bar{w}<
(fg)^n(g(\bar{w}))\,. \,$ Condition~$iii)$ follows because $\bar{w} <
fg^2(\bar{w})$ implies $fg^2(\bar{w}) <
fg^2(fg^2(\bar{w}))=(fg^2)^2(\bar{w})$, and $fg^2(\bar{w}) <
g(\bar{w})$ implies $(fg)^2(g(\bar{w}))= fg (fg^2(\bar{w})) < fg
(g(\bar{w}))=fg^2(\bar{w})$. $\hfill\square$

\vspace{0.5cm}


\beginpicture

\setcoordinatesystem units <1cm,1cm>


\putrule from 1.5 -2.5 to 6.5 -2.5 \putrule from 1.5 2.5 to 6.5 2.5
\putrule from 1.5 -2.5 to 1.5 2.5 \putrule from 6.5 -2.5 to 6.5 2.5

\plot 1.5 0 1.625 0.05296 1.75 0.109 1.875 0.16056 2 0.216 2.125
0.26056 2.25 0.3009 2.375 0.35296 2.5 0.4 /

\plot 2.5 0.4 2.625 0.45296 2.75 0.509 2.875 0.56056 3 0.616 3.125
0.66056 3.25 0.7009 3.375 0.75296 3.5 0.8 /

\plot 3.5 0.8 3.625 0.85296 3.75 0.909 3.875 0.96056 4 1.016 4.125
1.06056 4.25 1.1009 4.375 1.15296 4.5 1.2 /

\plot 4.5 1.2 4.625 1.25296 4.75 1.309 4.875 1.36056 5 1.416 5.125
1.46056 5.25 1.5009 5.375 1.55296 5.5 1.6 5.625 1.65296 /


\plot 6.5 -0.82519 6.375 -0.85296 6.125 -0.909 5.875 -0.96056 5.625
-1.016 5.375 -1.06056 5.125 -1.1009 4.875 -1.15296 4.625 -1.2 /

\plot 4.625 -1.2 4.375 -1.25296 4.125 -1.309 3.875 -1.36056 3.625
-1.416 3.375 -1.46056 3.125 -1.5009 2.875 -1.55296 2.625 -1.6 2.375
-1.65296 /


\setdots

\plot 1.5 -2.5 6.5 2.5 /

\putrule from 2.32 -1.68 to 5.68 -1.68 \putrule from 2.32 -1.68 to
2.32 1.68 \putrule from 2.32 1.68 to 5.68 1.68 \putrule from 5.68
-1.68 to 5.68 1.68

\put{Figure 1: The crossing} at 4 -3.5 \put{} at -4.2 0

\small


\put{$\bar{w}$} at 1.5 -2.8 \put{$g(\bar{w})$} at 6.5 -2.8
\put{$fg^2(\bar{w})$} at 4 -2.8  \put{$\bullet$} at 1.5 -2.5
\put{$\bullet$} at 6.5 -2.5 \put{$\bullet$} at 4   -2.5

\put{$fg$} at 4.7 -0.9 \put{$fg^2$} at 3.4 1.1

\endpicture


\vspace{0.5cm}

Note that if we let $w_0$ be the first element (w.r.t. $\leq^*$) of
$\Omega \,, \,$ then the stabilizer of $\, w_0 \,$ is $\,\preceq$-convex.
Indeed, if $id\prec g \prec f$, with $f(w_0)=w_0 \,, \,$ then $\, w_0 <^*
w_f \leq^* w_g$, and thus $g(w_0)=w_0\,. \,$ Actually, it is
not hard to see that the same argument shows the following.

\vsp \vsp

\begin{prop} \label{laclave} {\em Let $\Omega$ be a set endowed with a well-order
$\,\leq^* . \;$ If a group $G$ acts faithfully on $\Omega$ preserving a total order on it,
then there exists a left-ordering on $G$ for which the stabilizer $G_{\Omega_0}$
of any initial segment $\Omega_0 \,$ of $\, \Omega \,$ (w.r.t. $\leq^*$) is convex.
Moreover, if the action has no crossing, then this ordering is Conradian.}
\end{prop}

\vsp

\begin{ex} A very useful example of an action without crossings is
the action by left translations on the set of left-cosets of any
subgroup $H$ which is convex with respect to a $C$-ordering
$\,\preceq \,$ on $G$. Indeed, it is not hard to see that, due to
the convexity of $H$, the order $\, \preceq \,$ induces a total
order $\, \preceq_H \,$ on the set of left-cosets $G / H$. Moreover,
$\, \preceq_H \,$ is $G-$invariant. Now suppose that $\, (f, g, uH,
vH, wH) \,$ is a crossing for the action. Since $w_1H\prec_{H} w_2
H$ implies $\,w_1\prec w_2 \,, \,$ for all $w_1,\, w_2$ in $G$, we
have that $(f,g,u,v,w)$ is actually a crossing for the action by
left translations of $G$ on itself. Nevertheless, this contradicts
Theorem \ref{thm-navas}. \label{convex-cosets}
\end{ex}

\vsp

The following is an application of the preceding example. For
the statement, we will say that a subgroup $H$ of a group $G$ is
\textit{$C$-relatively convex} if there exists a $C$-ordering
on $G$ for which $H$ is convex.

\vsp

\begin{lem} {\em For every $C$-orderable group, the intersection of any
family of $C$-relatively convex subgroups is $C$-relatively convex.}
\label{int-conv}
\end{lem}

\noindent{\em Proof.} We consider the action of $G$ by left multiplications
on each coordinate of the set $\Omega=\prod_\alpha G / H_\alpha \,, \,$
where $(G / H_\alpha,\;\preceq_{H_\alpha})$ is the ($G-$invariant
ordered) set of left-cosets of the $C$-relatively convex subgroup
$H_\alpha$. Putting the (left) lexicographic order on $\Omega$
and using Example \ref{convex-cosets}, it is not hard to see
that this action has no crossing. Moreover, since $\{id\}$
is $C$-convex, the action is faithful.

\vsp

Now consider an arbitrary family $\Omega_0\subset\{H_\alpha\}_\alpha$
of $C$-relatively convex subgroups of $G$, and let $\leq^*$ be a well-order
on $\Omega$ for which $\Omega_0$ is an initial segment. For the induced
ordering $\preceq$ on
$G$, it follows from Proposition \ref{laclave} that the stabilizer
$G_{\Omega_0} = \bigcap_{H\in\, \Omega_0} H$ is $\preceq$-convex.
Moreover, Lemma \ref{generalization} implies that $\preceq$
is a $C$-ordering, thus concluding the proof. $\hfill \square$
\newline

We close this section with a simple lemma that we will need later
and which may be left as an exercise to the reader (see also
\cite[Lemma 5.2.3]{kopytov}).

\vsp

\begin{lem} {\em Let $G$ be a torsion-free Abelian group. Then $G$
admits only finitely many $C$-orderings if and only if $G$ has rank 1}.
\label{abel-1}
\end{lem}


\section{Proof of Theorem B}
\label{B}

\subsection{On groups with finitely many $C$-orderings}
\label{sec}

\hspace{0.45cm} Let $G$ be a $C$-orderable group admitting only
finitely many $C$-orderings. Obviously, each of these orderings must
be isolated in $\mathcal{CO}(G)$. We claim that, in general, if
$\, \preceq \,$ is an isolated $C$-ordering, then the series of
$\, \preceq$-convex subgroups
$$\{id\} =  G^{id} \subset \ldots \subset G_g\lhd G^g \subset \ldots \subset G$$
is finite. Indeed, let $\{f_1, \ldots, f_n\}\subset G$ be a set of
$\, \preceq$-positive elements such that $V_{f_1,\ldots,f_n}$ consists
only of $\, \preceq . \,$ If the series above is infinite, then there
exists a $g\in G$ so that no $f_i$ belongs to $G^g\setminus G_g$.
This implies that the flipped ordering $\,\preceq_g \,$ is Conradian and
different from $\,\preceq . \,$ However, every $f_i$ is still $\, \preceq_g$-positive
({\em c.f.}, Example \ref{exflip}), which is impossible
because $V_{f_1,\ldots,f_n} = \{\preceq\}$.

\vsp

Next let
$$\{id\} = G_0 \lhd G_1 \lhd \ldots \lhd G_n = G$$
be the series of $\, \preceq$-convex subgroups of $G.\,$ According to the Conrad Theorem,
every quotient $G_i/G_{i-1}$ embeds into $\mathbb{R}$,
and thus it is Abelian. Since every ordering on such a quotient
can be extended to an ordering on $G$ (similarly as in
Example \ref{exflip}), the Abelian quotient $G_i/G_{i-1}$
has only a finite number of orderings. It now follows from Lemma \ref{abel-1}
that it must have rank 1. Therefore, the series above is rational.

We now show that this series is unique. Suppose
$$\{id\}=H_0\lhd H_1\lhd \ldots \lhd H_k =G$$
is another rational series. Since $H_{k-1}$ is $C$-relatively convex,
we conclude that
$\, N=G_{n-1}\cap H_{k-1} \,$ is $\,C-$relatively convex by Lemma \ref{int-conv}. Now
$\,G/N \,$ is torsion-free Abelian and has only a finite number
of orderings, thus it has rank 1. Since convex groups are {\em isolated},
$H_{k-1}$ and $G_{n-1}$ have the property that $x^r \in G_{n-1}$
(resp. $x^r \in H_{k-1}$) implies $x \in G_{n-1}$ (resp. $x \in H_{k-1}$).
This clearly yields $H_{k-1} = G_{n-1}$. Repeating this argument several
times, we conclude the uniqueness of the rational series,
which is hence normal.

\vsp

Now we claim that no quotient $\, G_{i+2}/G_i \,$ is Abelian. If
not, $G_{i+2}/G_i$ would be a rank-2 Abelian group, and so
an infinite number of orderings could be defined on it. But
since every ordering on this quotient can be extended to a
$C$-ordering on $G$, this would lead to a contradiction.


\subsection{On groups with a normal rational series}
In this subsction we prove the converse of Theorem B.
\medskip

Suppose that $\,G \,$ has a normal rational series
$$\{id\}=G_0 \lhd G_1\lhd ... \lhd G_n=G$$
such that no quotient $\,G_{i+2}/G_i \,$ is Abelian. Clearly, flipping
the orderings on the quotients $\,G_{i+1}/G_i \,$ we obtain at least $2^n$
many $C$-orderings on $G$. We claim that these are the only possible
$C$-orderings on $G$. Indeed, let $\, \preceq \,$ be a $C$-ordering on $G$,
and let $a \in G_1$ and $b \in G_2\setminus G_1$ be two non-commuting
elements. Denoting the Conrad homomorphism of the group
$\, \langle a,b\rangle \,$ endowed with the restriction of $\, \preceq \,$
by $\, \tau , \,$ we have $\tau(a) = \tau(bab^{-1})$. Since $G_1$ is rank-1
Abelian, we have $bab^{-1}=a^r$ for some rational number $r\not=1$. Thus
$\tau(a)=r\tau(a)$, which implies that $\tau(a)=0$. Therefore, $\,a <<
|b|\,, \,$ or in other words
$\, a^n\prec |b| \,$ for every $n \in \mathbb{Z}\,, \,$
where $|b| = max\{b^{-1},b\} \,. \,$ Since $G_2/G_1$ is rank-1, this
actually holds for every $b \neq id$ in $G_2 \setminus G_1$. Thus
$G_1$ is convex in $G_2$.

\vsp

Repeating the argument above, though now
with $G_{i+1}/G_i$ and $G_{i+2}/G_i$ instead of $G_1$ and $G_2$,
respectively, we see that the rational series we began with is none other
than the series given by the convex subgroups of $\,\preceq . \,$ Since
each $G_{i+1}/G_i$ is rank-1 Abelian, if we choose
$b_i \in G_{i+1} \setminus G_i$ for each $i=0,\ldots,n-1$, then
any $C$-ordering on $G$ is completely determined by the signs
of these elements. This shows that $G$ admits precisely
$2^n$ different $C$-orderings.


\section{Proof of Theorem A}

\hspace{0.45cm} Let $G$ be a group admitting a $C$-ordering $\,\preceq \,$ which is
isolated in the space of $C$-orderings. As we have seen at the beginning of
\S \ref{sec}, the series of $\,\preceq$-convex subgroups must be finite, say
$$\{id\} = G_0 \lhd G_1 \lhd \ldots \lhd G_n = G.$$
Proceeding just as in Example \ref{exflip}, any ordering on
$\,G_{i+1}/G_i \,$ may be extended (preserving the set of positive
elements outside of it) to a $C$-ordering on $G$. Hence, each quotient
must be rank-1 Abelian, so the series above is rational. We claim that
this series of $\,\preceq$-convex subgroups is unique (hence normal) and that no quotient
$\,G_{i+2}/G_i \,$ is Abelian. In fact, if the series has length 2, then it is normal.
Moreover, since no $C$-ordering on a rank-2 Abelian group is
isolated, we have that $G_2$ is non-Abelian. Then by Theorem B, the series is
unique. In the general case, we will use induction on the length of
the series. Suppose that every group having an isolated $C$-ordering
whose rational series of convex subgroups
$$\{id\} = H_0 \lhd H_1 \lhd \ldots \lhd H_{k}$$
has length $k< n$ admits a unique (hence normal) rational series
and that no quotient $H_{i+2}/H_i$ is Abelian. Let
$$\{id\}=G_0\lhd  \ldots \lhd G_{n-2} \lhd G_{n-1} \lhd G_n=G$$
be a rational series of length $n$ associated to some isolated
$C$-ordering $\,\preceq \,$ on $G$. Since $G_{n-1}$ is normal in $G$,
for every $g\in G$, the conjugate series
$$\{id\}=G^g_0\lhd  \ldots \lhd G^g_{n-2} \lhd G^g_{n-1} = G_{n-1}$$
is also a rational series for $G_{n-1}$. Since this series
is associated to a certain isolated $C$-ordering, namely the restriction
of $\, \preceq \,$ to $G_{n-1}$, we conclude that it is unique by
the induction hypothesis. Hence the series must coincide with the
original one, or in other words $G_{i}^g = G_{i}$.
Therefore, the series for $G$ is normal. Moreover,
every quotient $G_{i+2}/G_i$ is non-Abelian, because if not then $\preceq$
could be approximated by other $C$-orderings on $G$. Thus, by Theorem B,
the rational series for $G$ is unique, and $G$ admits only finitely many
$C$-orderings. This completes the proof of Theorem A.

\vsp\vsp\vsp

We conclude this section with a short discussion on the structure
of bi-orderable groups admitting finitely many $C$-orderings.
To begin with, let us note the following simple

\vsp

\begin{prop} {\em If a group $G$ has a bi-order that is isolated in the space of
$C$-orderings, then $G$ has finitely many $C$-orderings, each of which is bi-invariant.}
\end{prop}

\noindent {\em Proof.} The fact that $G$ admits only finitely many $C$-orderings
is direct consequence of Theorem A. Let $\,\preceq \,$ be a bi-ordering on $G$ which
is isolated in the space of $C$-orderings, and let $\, \preceq' \,$ be any other
$C$-ordering on $G$. According to the proof of Theorem B, the series of convex subgroups
$$\{id\} = G_0 \lhd G_1 \lhd \ldots \lhd G_n = G$$
is the same for both $\, \preceq \,$ and $\,\preceq'$. Moreover,
$\, \preceq' \,$ is obtained from $\,\preceq \,$ by flipping
the ordering on some of the quotients $G_{i+1} / G_i$. Since
$\, \preceq \,$ is bi-invariant, the set of $\,\preceq$-positive elements
in $G_{i+1} \setminus G_i$ is invariant under conjugacy. Since the flipping
procedure interchanges the sets of positive and negative elements (where
{\em negative} means non-trivial and non-positive), the
above remains true after flipping. More precisely,
the set of $\,\preceq'$-positive elements in $G_{i+1} \setminus G_i$
is invariant under conjugacy. Since this is true for every index $i$,
this shows that $\,\preceq'\,$ is bi-invariant. $\hfill\square$

\vspace{0.4cm}

In spite of the fact that the preceding proposition holds for general $n \geq 1, \,$
it only applies for $n = 1,2$.

\vsp

\begin{prop} {\em If a bi-orderable group has finitely many
$C$-orderings, then the number of $C$-orderings is two or four.}
\end{prop}

\noindent{\em Proof.} Let $G$ be a bi-orderable group having $2^n$
$\,C-$orderings for some $n \geq 3$, and let
$$\{id\} = G_0 \lhd G_1 \lhd G_2 \lhd G_3 \unlhd \ldots \unlhd G_n=G$$
be the series of convex subgroups. It is easy to see that given $a
\in G_1$, there exist $b\in G_2\setminus G_1$ and $c\in G_3\setminus
G_2$ such that $bab^{-1}= a^r$ and $cb^qc^{-1}=b^{p}w$, where $w\in
G_1$, and $r, p/q$ are positive rational numbers such that $p/q\not=
1$. Let $t \in \mathbb{Q}$ be positive and such that $cac^{-1}=a^t$.
We have
$$a^{r^p} = b^p w aw^{-1} b^{-p}= (c b^q c^{-1}) a (c b^{-q}
c^{-1})= cb^q a^{1/t} b^{-q} c^{-1}  = ca^{r^q/t}c^{-1} = a^{r^q}.$$
Since $0 < r \neq 1$, we have that $p=q$, which contradicts
the fact that $1\not= p/q$. $\hfill\square$


\section{Some non-trivial examples}

\subsection{The Baumslag-Solitar group}
\label{ejem}

\hspace{0.45cm} In this subsection, we consider the Baumslag-Solitar
group $\,B(1,2) = \langle a, b \mid b a b^{-1} = a^2 \rangle$. We
let $\langle\langle \,a \rangle\rangle$ denote the largest rank-1
Abelian subgroup cointaining $a$. According to \cite[\S
5.3]{kopytov}, this group admits four bi-orderings which are
obtained from the series
$$\{id\} \esp \lhd \esp \langle\langle \,a \rangle\rangle = a^{\mathbb{Z}[\frac12]} \esp
\lhd \esp \langle a, b \rangle.$$
Here $\,\mathbb{Z}[\frac{1}{2}]\,$ denotes the set of dyadic rational numbers,
that is,
$$\mathbb{Z}\left[\frac{1}{2}\right]=\left\{\;\frac{m}{2^k}\;
\mid m \in \Z,\;k\in \N \right\}.$$

Below we give a self-contained proof of the fact that, actually, any
$C$-ordering on $B(1,2)$ coincides with one of these bi-orderings.

\vspace{0.2cm}

\begin{prop} {\em Baumslag-Solitar's group $B(1,2)$ admits only four $C$-orderings.}
\label{solo-cuatro}
\end{prop}

\noindent {\em Proof.} Let $\,\preceq \,$ be a $C$-ordering on
$B(1,2)$. We will determine all $\,\preceq$-convex subgroups of
$B(1,2)$. Without lost of generality, we may assume that $\, b \succ
1 \,$. Otherwise, we could change $\, \preceq \,$ by its `opposite'
ordering $\, \overline{\preceq} \,$ ({\em i.e.}, the one whose
positive elements are the inverses of the positive elements of $\,
\preceq \,$; compare Example \ref{exflip}) which has the same convex
subgroups.

\vsp

First we claim that $a << b$ ({\em i.e.}, $a^n\prec b$ for all
$n\in\mathbb{Z}$). Indeed, if we let $\tau$ be the (unique up to
multiplication by a positive real number) Conrad homomorphism, then
we have \esp $\tau(a)=\tau(b a b^{-1})=\tau(a^2)=2\tau(a)$ \esp
which implies that $\tau(a)=0$. Hence $\tau(a^n) = 0$ for all $n \in
\mathbb{Z}$. Since $\tau$ is non-trivial and non-decreasing, we must
have $\tau(b) > 0$. Again, from the fact that $\tau$ is
non-decreasing, we conclude that $a^n\prec b$ as asserted.

\vsp

Next letting $G_g \subset G^g$ be the convex jump associated to $g
\in B(1,2)$, by property (3) in the Introduction we have $a \in
G_{b} \neq B(1,2)$. It follows that an arbitrary $h = a^{n_1}
b^{m_1} \cdots a^{n_i} b^{m_i}$ is contained in $ker(\tau)$ if and
only if $\sum_k m_k = 0$. It is easy to check that $\sum_k m_k = 0$
holds if and only if $h \in a^{\mathbb{Z}[\frac{1}{2}]}$. This shows
that the sequence of $\preceq$-convex subgroups of $B(1,2)$ is
$$\{id\} = G_{a} \subset G^{a} = a^{\mathbb{Z}[\frac{1}{2}]}
= G_{b} \subset G^{b} = B(1,2).$$ Since both $G_{b}$ and $G^{b} /
G_{b}$ are torsion-free rank-1 Abelian, we have that $B(1,2)$ admits
only four $C$-orderings. $\hfill\square$

\vspace{0.45cm}

We point out that $B(1,2)$ admits infinitely many left-orderings.
Indeed, let $\xi\!: B(1,2)\to \mathrm{Homeo}_+ (\mathbb{R})$ be the
isomorphic imbedding given by $a \!: x \mapsto x+1$ and $b \!: x
\mapsto 2x$. We associate to each irrational number $\varepsilon$ a
left-ordering $\preceq_{\varepsilon}$ on $B(1,2)$ whose set of
positive elements is defined by $\{g \!\in\! B(1,2) \mid \esp \xi(g)
(\varepsilon)
> \varepsilon \}$. (These orderings were introduced by Smirnov in
\cite{smirnov}.) When $\varepsilon$ is rational, the preceding set
defines only a partial ordering. However, in this case the
stabilizer of the point $\varepsilon$ is isomorphic to $\mathbb{Z}$,
and hence this partial ordering may be completed to two total
left-orderings $\preceq_{\varepsilon}^+$ and
$\preceq_{\varepsilon}^{-}$. Here $\preceq_{\varepsilon}^+$ (resp.
$\preceq_{\varepsilon}^-$) corresponds to the limit of
$\preceq_{\varepsilon_n}$ for any sequence of irrational numbers
converging to $\varepsilon$ by the right (resp. left).

\vsp

Remark that the opposite orderings ({\em i.e.,} those of the form
$\overline{\preceq}_{\varepsilon}$) can be obtained the same way as
above though now starting with the embedding $a \!: x \mapsto x-1$,
$b \!: x \mapsto 2x$ (and changing $\varepsilon$ by $-\varepsilon$).
Moreover, as $\varepsilon$ tends to $-\infty$ or $+\infty$, the
associate orderings converge to bi-orderings. This corroborates a
result by Navas (see \cite[Proposition 4.1]{navas}) according to
which no $C$-ordering is isolated in the space of left-orderings of
a group having infinitely many left-orderings.

\vsp

We next give a complete description of the space of left-orderings
of $B(1,2)$.

\vsp\vsp

\begin{thm}\label{laprop} {\em Besides the four bi-orderings previously described, the
space of left-orderings of $B(1,2)$ is made up of those of the form $\preceq_{\varepsilon}$
for $\varepsilon \notin \mathbb{Q}$, those of the form $\preceq_{\varepsilon}^+$ and
$\preceq_{\varepsilon}^-$ for $\varepsilon \in \mathbb{Q}$, and their opposites.
In particular, every non Conradian ordering on $B(1,2)$ can be realized as an
induced ordering coming from an affine action of $B(1,2)$ on the real line.}
\end{thm}

\vsp

To prove this theorem, we will use the ideas involved in the
following well-known orderability criterion (see \cite[Theorem
6.8]{ghys}, \cite[\S 2.2.3]{book}, or \cite{navas} for further
details).

\vsp

\begin{prop} {\em For a countable infinite group $G$, the following two properties are
equivalent:\\

\noindent -- $G$ is left-orderable,\\

\noindent -- $G$ acts faithfully on the real line by orientation
preserving homeomorphisms.}
\end{prop}

\noindent\textit{Sketch of proof. } The fact that a group of
orientation preserving homeomorphisms of the real line is
left-orderable is a direct consequence of Theorem \ref{Cohn}.

\vsp

For the converse, we construct what is called \textit{the dynamical
realization of a left-ordering}. Let $\preceq$ be a left-ordering on
$G$. Fix an enumeration $(g_i)_{i \geq 0}$ of $G$, and let
$t(g_0)=0$. We shall define an order-preserving map $t: G \to \R$ by
induction. Suppose that $t(g_0), t(g_1),\ldots,t(g_i)$ have been
already defined. Then if $g_{i+1}$ is greater (resp. smaller) than
all $g_0,\ldots, g_i$, we define $t(g_{i+1})= max\{t(g_0),\ldots,
t(g_i)\}+1$ (resp. $min\{t(g_0),\ldots, t(g_i)\}-1$). If $g_{i+1}$
is neither greater nor smaller than all $g_0,\ldots,g_i$, then there
are $g_n,g_m\in\{g_0,\ldots , g_i \}$ such that $g_n\prec
g_{i+1}\prec g_m$ and no $g_j$ is between $g_n,g_m$ for $0\leq j\leq
i$. Then we put $t(g_{i+1})=(t(g_n)+t(g_m))/2$.

\vsp

Note that $G$ acts naturally on $t(G)$ by $g(t(g_i)) = t(gg_i)$. It
is not difficult to see that this action extends continuously to the
closure of $t(G)$. Finally, one can extend the action to the whole
real line by declaring the map $g$ to be affine on each
interval in the complement of $t(G)$. $\hfill\square$

\vsp

\begin{rem} As constructed above, the dynamical realization depends
not only on the left-ordering $\preceq$, but also on the enumeration
$(g_i)_{i\geq 0}$. Nevertheless, it is not hard to check that
dynamical realizations associated to different enumerations (but the
same ordering) are \textit{topologically conjugate}.\footnote{Two
actions $\phi_1\!: G \to \mathrm{Homeo}_+(\R)$ and $\phi_2\!:G \to
\mathrm{Homeo}_+(\R)$ are topologically conjugate if there exists
$\varphi \in \mathrm{Homeo}_+(\R)$ such that $\varphi\circ \phi_1(g)
= \phi_2(g) \circ \varphi$ for all $g \in G$.} Thus, up to
topological conjugacy, the dynamical realization depends only on the
ordering $\preceq$ of $G$.

\vsp

An important property of dynamical realizations is that they do not
admit global fixed points (\textit{i.e.,} no point is stabilized by
the whole group). Another important property is that $g\succ id$ if
and only if $g(t(id))> t(id)$, which allows us to recover the
left-ordering from the dynamical realization.
\end{rem}

\vsp

\noindent \textbf{Proof of Theorem \ref{laprop}. } Given a
left-ordering $\preceq$ on $B(1,2)$ we will consider its dynamical
realization. We have the following two cases:

\vspace{0.35cm}

\noindent \textbf{Case 1.} The element $a \in B(1,2)$ is cofinal
(that is, for every $g_1,\, g_2 \in B(1,2)$, there are $n_1,\, n_2
\in \mathbb{Z}$ such that $a^{n_1}\prec g_1$ and $a^{n_2}\succ
g_2$).

\vsp\vsp\vsp

For the next two claims, recall that for any measure $\mu$ on a
measurable space $X$ and any measurable function $f: X\to X$,
the {\em push-forward measure} $f_*(\mu)$ is
defined by $f_*(\mu)(A) = \mu(f^{-1}(A))$, where $A\subseteq X$
is a measurable subset. Note that $f_*(\mu)$
is trivial if and only if $\mu$ is trivial. Moreover,
one has $(fg)_*(\mu) = f_*(g_*(\mu))$ for all
measurable functions $f,g$.

\vsp\vsp\vsp\vsp

\noindent{\underbar{Claim 1.}} The subgroup $\langle\langle \,a
\rangle\rangle$ preserves a Radon measure $\nu$ ({\em i.e.,} a
measure which is finite on compact sets)  on the real line which is
unique up to scalar multiplication and has no atoms.

\vsp\vsp\vsp

Since $a$ is cofinal and $\langle\langle \,a \rangle\rangle$ is
rank-1 Abelian, its action on the real line is {\em free} (that is,
no point is fixed by any non-trivial element). By H\"older's theorem
(see \cite[Theorem 6.10]{ghys} or \cite[\S 2.2]{book}),
$\langle\langle \, a \rangle\rangle$ is semi-conjugated to a group
of translations. More precisely, there exists a non-decreasing,
continuous, surjective function $\varphi \!: \mathbb{R} \rightarrow
\mathbb{R}$ such that, to each $g \in \langle\langle \, a
\rangle\rangle$ one may associate a translation parameter $c_g$ so
that, for all $x \in \mathbb{R}$,
$$ \varphi(g(x))=\varphi(x)+ c_g.$$
Now since the Lebesgue measure $Leb$ on the real line is invariant
by translations, the {\em push-backward measure} $ \nu =
\varphi^*(Leb)$ is invariant by $\langle\langle \, a
\rangle\rangle$. (Here $\varphi^*(Leb)$ is defined by
$\varphi^*(Leb)(A) = Leb(\varphi(A))$.) Since $Leb$ is a Radon
measure without atoms, this is also the case for $\nu$. Finally, the
uniqueness of $\nu$ up to scalar multiple is an easy exercise (see
for instance \cite[Proposition 2.2.38]{book}).

\vsp\vsp\vsp\vsp

\noindent{\underbar{Claim 2.}} For some $\lambda\not= 1$, we have $b_*(\nu) = \lambda \nu$.

\vsp\vsp\vsp

Since $\langle\langle \,a \rangle\rangle \lhd B(1,2)$, for any $a'
\in \langle\langle \,a \rangle\rangle$ and all measurable $A \subset
\mathbb{R}$ we must have $b_*(\nu)(a' (A)) = \nu(b^{-1}a'
(A))=\nu(\bar{a}(b^{-1}(A)))=\nu(b^{-1}(A))=b_*(\nu)((A))$ for some
$\bar{a} \in \langle\langle \,a\rangle\rangle$. (Actually, $a' =
\bar{a}^2$.) Thus $b_*(\nu)$ is a measure that is invariant by
$\langle\langle \,a \rangle\rangle$. The uniqueness of the
$\langle\langle \,a \rangle\rangle$-invariant measure up to scalar
factor yields $b_*(\nu)=\lambda \nu$ for some $\lambda > 0$. Assume
for a contradiction that $\lambda$ equals 1. Then the whole group
$B(1,2)$ preserves $\nu$. Thus there is a {\em translation number
homomorphism} $\tau_{\nu} \!: B(1,2)\to \mathbb{R}$ defined by
$$ \tau_{\nu}(g)= \left\{ \begin{array}{c c} \nu([x,g(x)]) &
\text{ if $g(x) \geq x$, }  \\ -\nu([g(x), x]) & \text{ if $g(x) <
x. $}
\end{array} \right .$$ (one easily checks that this definition does
not depend on $x\in \R$). The kernel of $\tau_{\nu}$ must contain
the commutator subgroup of $B(1,2)$; since $a \in B(1,2)^\prime$,
this yields $\tau_{\nu}(a)=0$. Nevertheless, this is impossible,
since --as is easy to see- the kernel of $\tau_{\nu}$ coincides with
the set of elements having fixed points on the real line (see for
instance \cite[\S 2.2.5]{book}).

\vsp\vsp\vsp

By Claims 1 and 2, for each $g \in B(1,2)$ we have $g_*(\nu) = \lambda_g (\nu)$
for some $\lambda_g > 0$. Moreover, $\lambda_a = 1$ and $\lambda_b = \lambda$.
\vsp

Now, for $x\in \R$, let $F(x) = sgn(x - t(id)) \cdot
\nu([t(id),x])$. (Note that $F(t(id)) = 0$.) Semi-conjugating the
dynamical realization by $F$ yields a faithful representation $A \!:
B(1,2)\to \mathrm{Homeo}_+(\R)$ of $B(1,2)$ in the group of
(orientation-preserving) affine homeomorphisms of the real line.
More precisely, for all $g \in B(1,2)$ and all $x \in \mathbb{R}$ we
have
$$ F (g(x)) = A_g (F(x))$$
where the affine map $A_g$ is given by \esp
$$A_g(x) \esp = \esp \frac{1}{\lambda_g}
x + \frac{1}{\lambda_g} \nu([t(g^{-1}),t(id)])$$
(here we use the convention $\nu([x,y]) = - \nu ([y,x])$ for
$x > y$). For instance, if $x > t(id)$ and $g \in B(1,2)$
are such that $g(x) > t(id)$, then
\begin{eqnarray*}
F (g(x))
&=& \frac{1}{\lambda_g} F(x) + \frac{1}{\lambda_g} \nu([t(g^{-1}),t(id)]).
\end{eqnarray*}

The action $A$ induces a (perhaps partial) left-ordering
$\preceq_A$, namely $f \succ_A id$ if and only if $A_g (0) > 0$.
Clearly, if the orbit under $A$ of $0$ is free (that is, for every
non-trivial element $g\in B(1,2)$, we have $A_g(0)\not=0$), then
$\preceq_A$ is total and coincides with $\preceq$ (our original
ordering).

\vsp

Now assume that the orbit of $0$ is not free. (This may arise for
example when $\, t(id) \,$ does not belong to the support of $\nu$).
In this case, the stabilizer of $0$ under the action $A$ is
isomorphic to $\Z$. Therefore, $\preceq$ coincides with either
$\preceq_A^+$ or $\preceq_A^-$ (the definition of $\preceq_A^{\pm}$
is similar to the definition of $\preceq_{\varepsilon}^{\pm}$ given
above).

\vsp

Due to the discussion above, we need to determine all possible
embeddings of $B(1,2)$ into the affine group.

\vsp\vsp

\begin{lem} {\em Every faithful representation of $B(1,2)$
in the affine group is given by}
$$ a\sim \left( \begin{array}{c c} 1 & \alpha \\ 0 &1 \end{array}
\right), \;\;\; b\sim \left( \begin{array}{c c} 2&\beta  \\ 0&1
\end{array} \right)$$
for some $\alpha \not=0$ and $\beta \in \mathbb{R}$.
\end{lem}

\noindent \textit{Proof.} One easily checks that a correspondence as above
induces a faithful representation. Conversely, let
$$a\sim \left( \begin{array}{c c} s & \alpha \\ 0 &1 \end{array}
\right),\;\;\; b\sim\left( \begin{array}{c c} t & \beta \\ 0 &1
\end{array} \right)$$
be a representation. Then the following equality must hold:
$$a^2\sim \left( \begin{array}{c c} s^2 & s\alpha + \alpha \\ 0 &1 \end{array}
\right)= \left( \begin{array}{c c} s & \alpha t - s\beta+\beta \\ 0 &
1 \end{array} \right)\sim bab^{-1}.$$
Thus $s = 1$, $t = 2$, and
since the representation is faithful, $\alpha\not=0$. $\hfill \square$

\vspace{0.35cm}

Let $\alpha,\beta$ be such that $A_a(x)=x+\alpha$ and
$A_b(x)=2x+\beta$. We claim that if the stabilizer of $0$ under $A$
is trivial --which implies in particular that $\beta \!\neq\! 0$-- ,
then $\preceq_A$ (and hence $\preceq$) coincides with
$\preceq_{\varepsilon}$ if $\alpha > 0$ (resp.
$\overline{\preceq}_{\varepsilon}$ if $\alpha < 0$), where
$\varepsilon = \beta / \alpha $. Indeed, if $\alpha > 0$, then for
each $g = b^n a^r \in B(1,2)$ we have $A_g(0) = 2^n r\alpha + (2^n -
1) \beta$. Hence $A_g (0) > 0$ holds if and only if
$$2^n \beta/\alpha +2^n r >\beta/\alpha.$$
Letting $\varepsilon = \beta/\alpha$, one easily checks that the
preceding inequality is equivalent to $g\succ_\varepsilon id$. The
claim now follows.

\vsp

In the case the stabilizer of $0$ under $A$ is isomorphic to
$\mathbb{Z}$, similar arguments to those given above show that
$\,\preceq \,$ coincides with either $\,\preceq_{\varepsilon}^+$, or
$\, \preceq_{\varepsilon}^{-}$, or $\,
\overline{\preceq}_{\varepsilon}^+ \,$, or $\,
\overline{\preceq}_{\varepsilon}^{-} \,$, where $\varepsilon$ again
equals $ \beta / \alpha$.

\vspace{0.43cm}

\noindent\textbf{Case 2.} The element $a\in B(1,2)$ is not cofinal.

\vsp

In this case, for the dynamical realization of $\,\preceq \,, \,$ the set
of fixed points of $a$, denoted $Fix(a)$, is non-empty. We claim that
$b(Fix(a))=Fix(a)$. Indeed, given $x\in Fix(a)$, we have
$$a(b(x))=ab(x)=a^{-1}ba(x)=a^{-1}(b(x))\, .$$
Hence $a^2(b(x))=b(x)$, which implies that $a(b(x)) = b(x)$ as asserted.
Observe that since there is no global fixed point
for the dynamical realization, we must have $b(x)\not=x \,, \,$ for all $x
\in Fix(a) \,.$

\vsp

Now suppose that $b\succ id$ (otherwise, we may consider the
opposite ordering), and let $x_1 = inf \{x \!\in\! Fix(a) \mid x
> t(id)\}$. We claim that $b(x_1)>x_1$. Suppose not. Then $b(x_1) <
x_1$, but since $b(t(id)) = t(b) > t(id)$, we also have $b(x_1) >
t(id)$. Therefore, $b(x_1)$ is a fixed point of $a$ situated in
$(t(id),x_1)$, which is a contradiction.

\vsp

We now claim that $t(b)>x_1$. Indeed, if not, then we would have
$b(t(id)) = t(b) < x_1$. (Note that $t(b)$ cannot be equal to $x_1$,
since $x_1$ is fixed by $a$, but $B(1,2)$ acts freely on
$t(B(1,2))$.) Since $b(x_1) > x_1$, this would yield $b^{-1} (x_1)\!
\in (t(id),x_1)$. However, since $b^{-1}(x_1)$ belongs to $Fix(a)$,
this contradicts the definition of $x_1$.

\vsp

We next claim that $b(x_{-1})\geq x_1$, where $x_{-1} = \sup \{x
\!\in\! Fix(a) \mid x < t(id)\}$. Indeed, since $b(x_{-1})$ is a
fixed point of $a$, to show this it is enough to show that
$b(x_{-1})> x_{-1}$. This can be easily checked using similar
arguments to those above.

\vsp

We finally claim that $\langle\langle \,a \rangle\rangle$ is a
convex subgroup. First note that, by the definition of the dynamical
realization, for every $g\in B(1,2)$ we have $\,t(g)=g(t(id))$.
Then, it follows that for every $g\in \langle\langle \, a
\rangle\rangle$, $t(g) \!\in (x_{-1},x_1)$. Now let $m \in
\mathbb{Z}$ and $r,s$ in $\mathbb{Z} [\frac{1}{2}]$ be such that $id
\prec b^m a^r \prec a^s$. Then we have $t(id) < b^m (t(a^r)) <
t(a^s) < x_1$. Since $b(x_{-1}) \geq x_1$, this easily yields $m =
0$, that is, $b^m a^r \in \langle\langle \,a \rangle\rangle$.

\vsp

We have thus proved that $\langle\langle \,a \rangle\rangle$ is a
convex (normal) subgroup of $B(1,2)$. Since the quotient $B(1,2) /
\langle\langle \,a \rangle\rangle$ is isomorphic to $\mathbb{Z}$, an
almost direct application of the characterization (4) in the
Introduction shows that the ordering $\, \preceq \,$ is Conradian.
This concludes the proof of Theorem \ref{laprop}. $\hfill\square$

\vsp

\begin{rem} It follows from Theorem \ref{laprop} that no ordering is
isolated in $\mathcal{LO}(B(1,2))$. Thus this space is homeomorphic
to the Cantor set.  Moreover, the natural conjugacy action of
$B(1,2)$ on $\mathcal{LO}(B(1,2)$ is `almost' transitive. More
precisely, for any irrational $\varepsilon \,, \,$ the orbit of
$\, \preceq_\varepsilon \,$ under $B(1,2)$ is dense in the subspace
$V_a$ formed by the orderings for which $a$ is positive. This
easily follows from the fact that, for all $g\in B(1,2)$, we have
$g(\preceq_\varepsilon) = \esp \preceq_{g^{-1}(\varepsilon)}\!\!.$ The
complementary subspace $V_{a^{-1}}$ of $\mathcal{LO}(B(1,2))$ is
densely covered by the orbit of $\;\overline{\preceq}_\varepsilon.$
\end{rem}

\begin{rem} The above method of proof also gives a complete classification --up to
topological semiconjugacy-- of all actions of $B(1,2)$ by orientation-preserving
homeomorphisms of the real line (compare \cite{wilk,Na-sol}). In particular, all
these actions come from left-orderings on the group (compare with the comment
before Question 2.3 in \cite{navas}).
\end{rem}


\subsection{Examples of groups with $2^n$ Conradian orderings but
infinitely many left-orderings}
\label{ejem-n}

\hspace{0.45cm} The classification of groups having finitely many
left-orderings was obtained by Tararin and appears in \cite[\S
 5.2]{kopytov}. An example of a group having precisely $2^n$ orders
is \esp $T_n = \Z^n$ endowed with the product rule
$$(b_n, \ldots , b_1) \cdot (b'_n, \ldots, b'_1)=(b_n+b'_n, \, (-1)^{b'_n}
b_{n-1} + b'_{n-1}\, , \ldots , \,(-1)^{b'_2} b_1+b'_1).$$
A presentation for $T_n$ is
$$T_n\cong \langle a_n,\ldots,a_1 \mid R_n\rangle,$$
where the set of relations $R_n$ is
$$\;a_{i+1}a_i a_{i+1}^{-1} = a_i^{-1} \quad \mbox{if}
\quad i < n, \qquad \mbox{and }
\qquad a_ia_j=a_ja_i \quad \mbox{if} \quad |i - j| \geq 2.$$

\vsp

A very simple dynamical argument shows that if a group has finitely
many left-orderings, then each of these orderings is Conradian
\cite[Lemma 3.45]{navas}. However, it is natural to ask whether for
each $n \geq 2$ there are groups having precisely $2^n$ Conradian
orderings but infinitely many left-orderings. As we have seen in the
preceding section, for $n=2$ this is the case of the Bumslag-Solitar
group $B(1,2)$. This holds (with a very similar proof) for many
other metabelian left-orderable groups, as for example all
Baumslag-Solitar's groups $B(1,\ell)$ for $\ell \geq 2$. It turns
out that, in order to construct examples for higher $n$ and having
$B(1,\ell)$ as a quotient by a normal convex subgroup, we need to
choose an odd integer $\ell$. As a concrete example, for $n \geq 3$ we
endow \esp $C_n = \Z \times \Z[\frac{1}{3}] \times \Z^{n}$ \esp with
the multiplication
\begin{multline*}
\Big( c, \,\frac{m}{3^{k}}\,, a_n, \ldots, a_1\Big) \cdot
\Big(c^\prime, \,\frac{m^\prime}{3^{k^\prime}}\,, a_n^\prime, \ldots, a_1^\prime \Big) =\\
= \Big( c + c^\prime, \;3^{c}\frac{m^\prime}{3^{k^\prime}}+\frac{m}{3^{k}}\;,
(-1)^{m} a^\prime_n+a_n \,, (-1)^{a_n} a^\prime_{n-1} + a_{n-1} \ldots,\,
(-1)^{a_2} a_1^\prime + a_1 \Big).
\end{multline*}
Note that the product rule is well defined because if \esp
$m/3^k = \bar{m}/3^{\bar{k}}$, \esp then \esp $(-1)^{m}=(-1)^{\bar{m}}$
\esp (it is here where we use the fact that \esp $\ell = 3$ \esp is odd).

The group $C_n$ admits the presentation
$$C_n\cong \langle c,b, a_n,\ldots ,a_1 \mid  cbc^{-1}=b^3\,,\;
ca_i=a_ic\,,\; ba_nb^{-1}= a_n^{-1}\,,\; ba_i=a_ib \esp \text{ if } \esp i\not=n
\,  ,\;R_n \rangle.$$
This group satisfies the hypotheses of Theorem B and has exactly
$2^{n+2}$ Conradian orderings. However, it has $B(1,3)$ as a
quotient by a normal convex subgroup. Since $B(1,3)$ admits
uncountably many left-orderings, the same is true for $C_n$.


\vsp\vsp \vspace{0.4cm}

\noindent{\bf Acknowledgments.} It is a very big pleasure to thank Andr\'es Navas for
introducing me to this beautiful theory and for many fruitful discussions on the subject.
It is also a pleasure to thank the anonymous referee for his patient consideration of
this paper (and an earlier version of it), as well as for his remarks, comments and
corrections. This work was partially funded by the PBCT-Conicyt Research
Network on Low Dimensional Dynamics.


\begin{small}

\vspace{0.1cm}


\vspace{0.37cm}

\noindent Crist\'obal Rivas\\

\vspace{-0.04cm}

\noindent Dep. de Matem\'aticas, Fac. de Ciencias, Univ. de Chile\\

\vspace{-0.04cm}

\noindent Las Palmeras 3425, \~Nu\~noa, Santiago, Chile\\

\vspace{-0.04cm}

\noindent Email address: cristobalrivas@u.uchile.cl

\end{small}



\begin{thebibliography}{Dillo 83}

\bibitem{botto} {\sc R. Botto-Mura \& A. Rhemtulla.} {\em Orderable groups.} Lecture Notes
in Pure and Applied Mathematics, Vol. {\bf 27}. Marcel Dekker, New York-Basel (1977).

\bibitem{wilk} {\sc L. Bursler \& A. Wilkinson.} Global rigidity of solvable group
actions on $\mathrm{S}^1$. {\em Geometry and Topology} {\bf 8} (2004), 877-924.

\bibitem{butts} {\sc R. Buttsworth.} A family of groups with a countable infinite
number of full orders. {\em Bull. Austr. Math. Soc.} {\bf 12} (1971), 97-104.

\bibitem{clay}{\sc A. Clay.} Isolated points in the space of left-orderings
of a group. Preprint (2008), arxiv:0812.2499.

\bibitem{conrad} {\sc P. Conrad.} Right-ordered groups. {\em Mich. Math. Journal}
{\bf 6} (1959), 267-275.

\bibitem{ghys} {\sc \'E. Ghys.} Groups acting on the circle.
\textit{L'Enseignement Math\'ematique} {\bf 47} (2001), 239-407.

\bibitem{holder} {\sc O. H\"older.} Die Axiome der Quantit\"at und die Lehre vom Mass.
{\em Verh. S\"achs. Ges. Wiss. Leipzig, Math. Phys.} {\bf 53} (1901), 1-64.

\bibitem{leslie} {\sc L. Jim\'enez.} {\em Din\'amica de grupos
ordenables.} Master thesis, Univ. de Chile (2007)

\bibitem{kopytov} {\sc V. Kopitov \& N. Medvedev.} {\em Right ordered groups.}
Siberian School of Algebra and Logic, Plenum Publ. Corp., New York (1996).

\bibitem{linnell} {\sc P. Linnell.} The topology on the space of left-orderings
of a group. Preprint (2006), arxiv:math/0607470.

\bibitem{witte} {\sc D. Morris-Witte.} Amenable groups that act on the line.
{\em Algebr. Geom. Topol.} {\bf 6} (2006), 2509-2518.

\bibitem{book} {\sc A. Navas.} {\em Groups of circle diffeomorphisms.}
Forthcoming book, arxiv:math/0607481v2. Spanish version published in
Ensaios Matem\'aticos, Braz. Math. Soc. (2007).

\bibitem{navas} {\sc A. Navas-Flores.} On the dynamics of left-orderable groups.
Preprint (2007), arxiv:0710.2466.

\bibitem{thompson} {\sc A. Navas \& C. Rivas.} Describing all bi-orderings
on Thompson's group F. To appear in {\em Groups, Geometry, and
Dynamics}, arxiv:0808.1688.

\bibitem{crossings} {\sc A. Navas \& C. Rivas.}, with an Appendix by {\sc A. Clay.}
A new characterization of Conrad's property for group orderings,
with applications. Preprint (2008), arxiv:0901.0880.

\bibitem{Na-sol} {\sc A. Navas.} Groupes r\'esolubles de diff\'eomorphismes de
l'intervalle, du cercle et de la droite. {\em Bull. of the Brazilian Math. Society}
{\bf 35} (2004), 13-50.

\bibitem{rr} {\sc A. Rhemtulla \& D. Rolfsen.} Local indicability in ordered groups:
braids and elementary amenable groups. {\em Proc. Amer. Math. Soc.} {\bf 130} (2002),
2569-2577.

\bibitem{sikora} {\sc A. Sikora.} Topology on the spaces of orderings of groups.
{\em Bull. London Math. Soc.} {\bf 36} (2004), 519-526.

\bibitem{smirnov} {\sc D. Smirnov.} Right orderable groups. {\em Algebra i Logika}
{\bf 5} (1966), 41-69.

\bibitem{zenkov} {\sc A. Zenkov.} On groups with an infinite set of
right orders. {\em Sibirsk. Mat. Zh.} {\bf 38} (1997), 90-92. English
translation: {\em Siberian Math. Journal} {\bf 38} (1997), 76-77.

\end{thebibliography}
\end{document}